\documentclass[12pt]{article}
\usepackage{graphicx}
\usepackage[sans]{dsfont}
\usepackage[english]{babel}
\usepackage{xfrac}
\usepackage{color}
\usepackage{amsfonts, amsthm, amsmath}

\usepackage[hypertexnames=false,
backref=page,
    pdftex,
    pdfpagemode=UseNone,
    breaklinks=true,
    extension=pdf,
    colorlinks=true,
    linkcolor=blue,
    citecolor=blue,
    urlcolor=blue,
]{hyperref}

\renewcommand{\Im}{\operatorname{Im}}

\newtheorem{theorem}{Theorem}

\title{Hermitian metrics on complex non-Kähler manifolds\thanks{
This survey expands the topic of a seminar given by the author at the Second International Conference on Differential Geometry in Fes in October 2024. The author would like to warmly thank the organizers for the kind invitation and the warm hospitality. Many thanks in particular to Professor Mohamed Abbassi.
This survey is based on joint collaborations with Simone Calamai, Mauricio Correa, Francesco Pediconi, Cristiano Spotti, Valentino Tosatti among others. The author would like to warmly thanks them, as well as Lucia Alessandrini, Giuseppe Barbaro, Giovanni Bazzoni, Eleonora Di Nezza, Filippo Fagioli, Elia Fusi, Alexis Garcia, Vincent Guedj, Chinh Lu, Marco Miceli, Antonio Otal, Maurizio Parton, Giovanni Placini, Tancredi Schettini Gherardini, Jonas Stelzig, Nicoletta Tardini, Adriano Tomassini, Luis Ugarte, Raquel Villacampa, Victor Vuletescu, Joshua Windare, Michela Zedda, for many stimulating discussions that have helped him grow and for their company on this journey. Thanks also to Jeff Streets and Fangyang Zheng for some useful comments and suggestions. 
The author was partially supported by project PRIN2022 ``Real and Complex Manifolds: Geometry and Holomorphic Dynamics'' (code 2022AP8HZ9) and by GNSAGA of INdAM.
}}
\author{Daniele Angella\thanks{Dipartimento di Matematica e Informatica, Università di Firenze, viale Morgagni 67/A, 50134 Firenze, Italy.
Email: \href{mailto:daniele.angella@unifi.it}{\texttt{daniele.angella (at) unifi.it}}
}}

\date{\today}

\begin{document}

\maketitle

\begin{center}{\itshape
Dedicated to Professor Paul Gauduchon, on his 80th birthday.\\
With deep gratitude for his contributions to mathematics and to the scientific community, and for his constant presence and inspiration throughout my professional journey.
}
\end{center}

\begin{abstract}
In this survey, we consider various analytic problems related to the geometry of the Chern connection on Hermitian manifolds, such as the existence of metrics with constant Chern-scalar curvature, generalizations of the K\"ahler-Einstein condition to the non-K\"ahler setting, and the convergence of the Chern-Ricci flow on compact complex surfaces.
\end{abstract}

\tableofcontents

\section{The holomorphic landscape}
We are interested in studying the geometry and topology of holomorphic manifolds. We prefer to use the term ``holomorphic manifold'' \cite{debartolomeis}, instead of the more common ``complex manifold'', since they behave as the objects in the category of holomorphic maps. More precisely, a {\em holomorphic manifold} of complex dimension $n$ is a topological space, with the Hausdorff and second-countable properties, such that any point has a neighbourhood homeomorphic to an open set of $\mathds C^n$, which is called ``local holomorphic chart'', such that the transition functions are biholomorphisms.
A {\em holomorphic map} between holomorphic manifolds is a function which is locally expressed by holomorphic functions.

\subsection{Almost-complex structures}
Given a local holomorphic chart $\varphi \colon U \stackrel{\simeq}{\to} \varphi(U) \subseteq \mathds C^n$ on $U\subseteq X$, with coordinates $\varphi=(z^1,\ldots,z^n)$, let $z^i = x^i + \sqrt{-1} y^i$ denote its real and imaginary parts. We have a well-defined endomorphism $J$ of the tangent bundle, locally given by
$$ J \sfrac{\partial}{\partial x^i} = \sfrac{\partial}{\partial y^i}, \qquad J \sfrac{\partial}{\partial y^i} = -\sfrac{\partial}{\partial x^i}. $$
This endomorphism satisfies $J^2=-\mathrm{id}$ and encodes the pointwise linear complex structure on the tangent spaces, varying smoothly with the point. An endomorphism $J\in \mathrm{End}(TX)$ satisfying $J^2=-\mathrm{id}$ is called an {\em almost-complex structure}.
Conversely, an almost-complex structure $J$ is naturally associated with local holomorphic coordinates if and only if the bundle of the $\sqrt{-1}$-eigenspaces, defined as $T^{1,0}_xX := \{ v-\sqrt{-1}Jv : v \in T_xX \} \subset T_xX \otimes_{\mathds R} \mathds C$, is involutive \cite{newlander-nirenberg}. This condition is equivalent to require that $[V^{1,0},W^{1,0}]^{0,1}=0$ for any vector fields $V$ and $W$, where $V=V^{1,0}+V^{0,1}$ denotes the decomposition with respect to $TX \otimes_{\mathds R}\mathds C = T^{1,0}X \oplus T^{0,1}X$ with $T^{0,1}X = \overline{T^{1,0}X}$. The tensor $N_J(V,W)=[V^{1,0},W^{1,0}]^{0,1}$ is called the {\em Nijenhuis tensor} of $J$. In dimension greater than four, there are no known examples of almost-complex manifolds that do not admit any (integrable) complex structure. For instance, it is an open problem whether the six-dimensional sphere $S^6$ admits an integrable complex structure, see \cite{agricola-bazzoni-goertsches-konstantis-rollenske} and the references therein. See also \cite{albanese-milivojevic} for a conjecture by Sullivan on the minimal sum of Betti numbers of a compact complex $n$-fold, with $n \geq 3$ being four.

\subsection{Complex cohomologies}
Under the action of an (almost-)complex structure $J$, the tangent bundle decomposition $TX \otimes_{\mathds R}\mathds C = T^{1,0}X \oplus T^{0,1}X$ induces a corresponding decomposition for its dual (the cotangent bundle) and its exterior powers. More precisely, the bundle of differential forms decomposes into bigraded components as $\wedge^k T^*X \otimes_{\mathds R}\mathds C = \bigoplus_{p+q=k} \wedge^{p,q}X$, where $\wedge^{p,q}X:=\wedge^p T^{*\,1,0}X \otimes \wedge^q T^{*\,0,1}X$. Denote by $\Omega^{p,q}(X)$ the space of the corresponding sections. Accordingly, the exterior differential $d \colon \Omega^{p,q}(X) \to \Omega^{p+q+1}T^*(X)\otimes_{\mathds R}\mathds C$ on $(p,q)$-differential forms decomposes into components. In general, the exterior differential has four components: $d \colon \Omega^{p,q}(X) \to \Omega^{p+2,q-1}(X) \oplus \Omega^{p+1,q}(X) \oplus \Omega^{p,q+1}(X) \oplus \Omega^{p-1,q+2}(X)$. When $J$ is integrable, only $\partial \colon \Omega^{p,q}(X) \to \Omega^{p+1,q}(X)$ and $\overline\partial \colon \Omega^{p,q}(X) \to \Omega^{p,q+1}(X)$ remain non-zero. The condition $d^2=0$ assures that $\partial^2=\partial\overline\partial+\overline\partial\partial=\overline\partial^2=0$. In other words, $(\Omega^{\bullet,\bullet}(X),\partial,\overline\partial)$ has the
structure of a double complex.
It follows that several cohomologies can be defined: the Dolbeault cohomology $H^{\bullet,\bullet}_{\overline\partial}(X)=\sfrac{\ker\overline\partial}{\mathrm{im}\,\overline\partial}$ is the cohomology of the sheaf $\mathcal O_X$ of germs of holomorphic functions; its conjugate $H^{\bullet,\bullet}_{\partial}(X)$ is defined similarly; the Bott-Chern cohomology $H^{\bullet,\bullet}_{BC}(X)= \sfrac{\ker d}{\mathrm{im}\,\partial\overline\partial}$; and its dual, the Aeppli cohomology $H^{\bullet,\bullet}_A(X)=\sfrac{\ker\partial\overline\partial}{\mathrm{im}\,\partial+\mathrm{im}\,\overline\partial}$.
For laziness, we refer to \cite{angella-LNM} for more details on the cohomology of almost-complex or holomorphic manifolds, although better references can be found elsewhere.
We refer to \cite{stelzig-JLMS, stelzig-AdvMath22, stelzig-AdvMath25}, as well as other works by Jonas Stelzig, for recent in-depth advances in the understanding of the double complexes of holomorphic manifolds. For a comprehensive review of recent results on non-K\"ahler Hodge theory and deformations of complex structures, we refer to \cite{popovici-book}.

\subsection{Hermitian metrics}
Thanks to the paracompactness property, we can take the (pullback of) the standard Hermitian inner product on each holomorphic chart and then glue them together using a partition of unity. This construction yields a {\em Hermitian metric} on $X$, defined as a family of Hermitian inner products $h_x$ on each complex vector space $(T_xX, J_x)$, varying smoothly with $x\in X$.
Via the linear isomorphisms $(T_xX, J_x) \simeq (T^{1,0}_xX, \sqrt{-1}) \simeq (T^{0,1}X, -\sqrt{-1})$, this induces a family of Hermitian inner products $h$ on each fibre $(T_xX,J_x)$. Decomposing this inner product into real and imaginary parts, say $\sfrac{1}{2}(g-\sqrt{-1}\omega)$, we obtain: a Riemannian metric $g$ on $X$, for which $J$ acts as isometry; a $J$-invariant $2$-form $\omega$ on $X$, satisfying $\omega(V,JV)>0$ pointwise for any nowhere-zero vector field $V$. Since $\omega=g(J\_,\_)$, we will refer to a Hermitian metric either as $g$ or $\omega$ without ambiguity.

\subsection{Complex algebraic geometry and K\"ahler geometry}
The first examples of {\em compact} holomorphic manifolds are the complex torus $\sfrac{\mathds C^n}{\mathds Z^n}$ and the {\em complex projective space} $\mathds CP^n = \sfrac{\mathds C^{n+1}\setminus 0}{\mathds C\setminus 0}$. Another fundamental construction is as follows. Consider a homogeneous polynomial in $n+1$ variables with no multiple roots. Its zero set in $\mathds CP^n$ defines a compact holomorphic manifold of complex dimension $n-1$. More generally, the zero set of a finite number of homogeneous polynomials defines a compact holomorphic manifold, provided that the polynomials are sufficiently generic to assure smoothness. Many geometric properties of such manifolds are indeed encoded in the ideal generated by its defining polynomials (see {\itshape e.g.} \cite{hubsch}), which explains the name {\em algebraic projective manifolds}. For example, a smooth hypersurface $X$ of degree $d$ in $\mathds CP^n$ has canonical bundle isomorphic to $\mathcal O_X(-n-1+d)$. In particular, for $d=n+1$, the hypersurface $X$ is a Calabi-Yau manifold, meaning it has trivial canonical bundle.

The examples $X$ discussed above share a common property: there exists a Hermitian metric $g$ that osculates to order $2$ the standard Hermitian inner product of $\mathds C^n$. Specifically, this means that for every point $x \in X$, there exist local holomorphic coordinates $(z^1, \ldots, z^n)$ on a neighbourhood $U$ such that $g = \sum_{i,j=1}^n (\delta_{ij}+o(z)) \, dz^i \otimes d\bar z^j$ at $x$.
Such metrics are known as {\em K\"ahler metrics}, see \cite{schouten-vandantzig, kahler}, since \cite{weil}.
Equivalently, a Hermitian metric is K\"ahler if and only if its associated $2$-form $\omega$ is symplectic, meaning that $d\omega=0$.
Indeed, on $\mathds CP^n$ with homogeneous coordinates $[z^0:z^1:\cdots:z^n]$, the Fubini-Study metric \cite{fubini, study}, defined as $\omega_{FS}= \sfrac{1}{4\pi}dd^c\log(\sum_i|z^i|^2)$, is a K\"ahler metric. This metric naturally induces a K\"ahler metric on any submanifolds.
More precisely, K\"ahler geometry is a sort of transcendental analogue \cite{demailly-analytic} of algebraic geometry. More precisely, it is well-known \cite{kodaira} that a compact K\"ahler manifold is an algebraic projective manifold if and only if its associated $2$-form $\omega$ defines a rational cohomology class $[\omega]\in H^2(X,\mathds Q)$.
K\"ahler geometry lies at the intersection of complex, Riemannian, and symplectic geometries. Each of these perspectives provides specific tools, and their compatibility allows for easily switching between these points of view. This interplay makes K\"ahler geometry nearly as powerful as algebraic geometry itself.

Since K\"ahler geometry ``represents a perfect synthesis of the Symplectic and the Holomorphic worlds'', it is meaningful to perform ``a sort of chemical analysis of symplectic and holomorphic contribution [...] in order to better understand the role of the different components of the theory'' \cite{debartolomeis-tomassini}
For this reason, we are interested in investigating the broad class of complex {\em non-K\"ahler} manifolds.

\section{Riemannian metrics on differentiable manifolds}
Let us first consider the Riemannian context, referring to \cite{besse} for further details and a comprehensive discussion. Let $g$ be a Riemannian metric on a differentiable manifold $M$. There exists a unique affine connection $D$ on $M$ that preserves $g$ ({\itshape i.e.}, $D g=0$) and is torsion-free ({\itshape i.e.}, $T=0$ where $T(V, W) := D_V W - D_W V - [V, W]$). This connection is known as the {\em Levi-Civita connection} and it is fully determined by the Koszul formula: $2g(D_VW,Z) = Vg(W, Z) + W g(Z, V) - Zg(V, W) + g([V, W], Z) - g([W, Z], V) + g([Z, V], W)$.
It defines a {\em curvature} tensor $R(V,W)= D_{[V,W]}-[D_V, D_W]$, which can be interpreted as a $\mathrm{End}(TX)$-valued $2$-form.
We can also interpret the curvature in the following way. Consider $\alpha$ to be a differential $k$-form with values in $TM$. Define the exterior derivative $d^D$ associated with $D$ by $(d^D\alpha)(V_0,\ldots,V_k)=\sum_i (-1)^i D_{V_i}(\alpha(V_0, \ldots, \hat{V_i}, \ldots, V_n)) + \sum_{i\neq j} (-1)^{i+j} \alpha([V_i,V_j], V_0, \ldots, \hat{V_i}, \ldots, \hat{X_j}, \ldots, X_n)$, where ``$\hat{X_i}$'' means skipped. Then $R=-(d^D)^2$ measure the failure of $d^D$ to be a genuine differential.
The curvature tensor satisfies several additional symmetries. It is clear that $R(V,W)=-R(W,V)$, and we also have the identity $g(R(V, W )Z, Y) = -g(R(V, W) Y, Z)$. Moreover, $R\wedge \mathrm{id}=0$ (known as the first, or algebraic, Bianchi identity) and $d^D R = 0$ (the second, or differential, Bianchi identity), see \cite{gauduchon-book} for details.
The curvature tensor contains exactly the information about the {\em sectional curvatures} $K(V,W)= \sfrac{g(R(V,W)V,W)}{(g(V,V)g(W,W)-g(V,W)^2)}$ of any plane $\sigma_p=\mathrm{span}\{V,W\} \subseteq T_pX$ at $p$, which correspond to the Gaussian curvature of the surface whose tangent plane at $p$ is $\sigma_p$.
Still, interesting information is encoded in the {\em Ricci curvature}, defined as $\mathrm{Ric}(V,W):=\mathrm{tr}(Z \mapsto R(V,Z)W)$.
Lastly, another important metric invariant is the function $s=\mathrm{tr}_gr$, called the {\em scalar curvature}.

We aim to find some {\em canonical metrics}, unique up to automorphisms, that would encode the topological or differentiable properties of the manifold.

\subsection{Space forms}
A natural choice for this is to look for Riemannian manifolds with `constant' curvature.
Asking for complete Riemannian metrics whose sectional curvature is constant (when varying the planes at $p$ and the point $p$, the latter being a consequence of the former in dimension $n\geq 3$, thanks to the differential Bianchi identity) leads to the notion of {\em space form}. There are only three simply-connected model space forms in any dimension $n$, up to homothety: the Euclidean metric on $\mathds R^n$ has flat curvature $R=0$; the unit sphere $S^n$ with the round metric has constant sectional curvature $+1$; and the hyperbolic space $H^n$ has constant sectional curvature $-1$.

\subsection{Einstein metrics}
Since the constant sectional curvature problem seems too rigid, one might consider weakening it to the condition `constant Ricci tensor'. This can be interpreted as the problem of searching for metrics $g$ such that $$\mathrm{Ric}_g=\lambda\cdot g$$ with $\lambda$ constant. Indeed, in dimension $n\geq 3$, even if $\lambda$ is a scalar function, it will result in a constant, given by a multiple of the scalar curvature, precisely $n\lambda=s$. (For surfaces, namely the case $n=2$, there is only one notion of curvature: the Riemann tensor, the Ricci tensor and the scalar curvature all coincide with the Gauss curvature.)
The problem is related to Einstein field equations for gravitation, which are given by $\mathrm{Ric} - \sfrac{1}{2}sg = T$, where $T$ is the energy-momentum tensor. Indeed, solutions in the vacuum ($T=0$) correspond to Ricci flat metrics.
The Einstein problem is a variational problem: more precisely,  Einstein metrics are characterized as the critical points for the functional $S(g)=\int_M s_g \, \mu_g$, (where $\mu_g$ denotes the Riemannian volume form,) which corresponds to the total scalar curvature, on the space of metrics with volume one. As Besse remarks \cite[page 6]{besse}, ``[d]espite the simplicity of the [Einstein condition] the reader should not imagine that examples are easy to find''. Indeed, ``[t]he author will be happy to stand you a meal in a starred restaurant in exchange of'' some new examples!
As far as now, the known examples of Einstein metrics are mainly constructed either by exploiting large symmetry group or through coupled equation.
In the homogeneous setting, Einstein metrics correspond to positive definite solutions of a system of quadratic equations. For instance, they exist on irreducible symmetric spaces, including compact semi-simple groups. A complete classification  of compact homogeneous Einstein spaces is available up to dimension $7$, and examples are known to exist in all dimensions up to $11$, with a counterexample appearing in dimension $12$. Recently, the Alekseevskii conjecture, stating that any homogeneous Einstein manifold of negative scalar curvature is diffeomorphic to $\mathds R^n$ and is a solvmanifold, was proven by \cite{bohm-lafuente}.
For more details, see the surveys \cite{besse, wang-1, wang-2, lauret, jablonski} and references therein.
The first non-homogeneous example of compact Einstein manifold with positive scalar curvature was constructed by Page \cite{page, page-pope} on the connected sum $\mathds CP^2 \sharp \mathds CP^2$ and later generalized by B\'erard-Bergery \cite{berardbergery}. By exploiting the theory of Riemannian submersion \cite{oneill} and the cohomogeneity one property, the Einstein equation can be reduced to a system of second-order ODEs.
Lastly, another powerful tool for producing Einstein metrics is the Aubin-Yau theorem \cite{yau, aubin} solving the Calabi conjecture. This theorem guarantees the existence of K\"ahler-Einstein metrics on compact complex K\"ahler manifolds $X$ with non-positive first Chern class $c_1(X)\leq 0$. Also the existence of K\"ahler-Einstein metrics on Fano manifolds (namely, with $c_1(X)>0$) has been recently understood \cite{chen-donaldson-sun-1, chen-donaldson-sun-2, chen-donaldson-sun-3, tian}.

Note that simply-connected manifolds with holomorphically trivial canonical bundle, known as {\em Calabi-Yau manifolds}, possess strong geometric properties that make them the most promising candidates for internal spaces in string theory compactifications \cite{candelas-horowitz-strominger-witten}. We refer to \cite{picard-Matrix} for a very recent survey on Calabi-Yau threefolds. Calabi-Yau manifolds, which admit K\"ahler metrics with holonomy in $\mathrm{SU}(n)$, play a crucial role in many string theory compactifications because they admit covariantly constant spinors, thereby preserving a portion of the original supersymmetry in the four-dimensional effective theory. Their geometric properties (Ricci-flatness, trivial canonical bundle) fit perfectly with the low-energy string/supergravity equations of motion and the requirement of having unbroken supersymmetry.
It is estimated that there are at least on the order of $10^{100}$ vacua matching the Standard Model gauge group and low-energy spectrum, see {\itshape e.g.} \cite{douglas}.
Therefore, machine learning techniques seem promising for learning topological and holomorphic invariants in large datasets, see {\itshape e.g.} \cite{he-LNM, alawadhi-angella-leonardo-schettinigherardini, hirst-schettinigherardini} and references therein.

We notice that Einstein metrics are stationary points of the normalized {\em Ricci flow}, which evolves a Riemannian metric according to $$\frac{\partial}{\partial t}g(t)=-2\mathrm{Ric}_{g(t)}+\frac{2}{n}r_{g(t)}\cdot g(t),$$ where $n=\dim M$ and $r_g=\sfrac{(\int_M s \mu_g)}{(\int_M \mu_g)}$ denotes the normalized average scalar curvature. It is a foundational result by Hamilton \cite{hamilton} that, on compact three-dimensional Riemannian manifolds with positive Ricci curvature, the volume normalized Ricci flow exists for all time and converges to an Einstein metric.
It is well-known that the Ricci flow preserves the K\"ahler condition, in which case it is referred to as the {\em K\"ahler-Ricci flow}. In his seminal work, Cao \cite{cao} reproved the Aubin-Yau theorem using the K\"ahler-Ricci flow, by showing that, for compact K\"ahler manifolds $X$ with $c_1(X) = 0$ (respectively, $c_1(X) < 0$), the (normalized) K\"ahler-Ricci flow starting at any K\"ahler metric $\omega_0$ smoothly converges to the unique K\"ahler-Einstein metrics in $[\omega_0]$ (respectively, in $-c_1(X)$).

\subsection{Yamabe problem}
Lastly, one may be interested in studying {\em constant scalar curvature metrics}, which is a weaker condition than Einstein. The problem, originally introduced by Yamabe \cite{yamabe}, was solved through the combined efforts of Trudinger \cite{trudinger}, Aubin \cite{aubin-2}, and Schoen \cite{schoen}, with Schoen utilizing the positive mass theorem of general relativity \cite{schoen-yau}. This effort culminated in proving that any Riemannian metric $g$ on a compact $n$-dimensional manifold can be conformally changed to a metric $u^{\sfrac{4}{(n-2)}}\cdot g$ with constant scalar curvature. The problem corresponds to finding a positive smooth solution to the semi-linear pde
$\frac{4(n-1)}{n-2}\Delta_g u-s_g u+c u^{\sfrac{(n+2)}{(n-2)}}=0$ where $\Delta_g$ denotes the Laplace-Beltrami operator of $g$ and $c$ is a constant.
In fact, the scalar curvature can be prescribed to be any possible function. More precisely, Kazdan and Warner \cite{kazdan-warner} proved that, on a compact manifold of dimension three or greater, any smooth function that attains a negative value somewhere is the scalar curvature of some Riemannian metric. Moreover, if the manifold admits metrics with strictly positive scalar curvature, then any smooth function is the scalar curvature of some Riemannian metric.

When the metric is K\"ahler, the scalar curvature has an interesting interpretation  as a moment map for an infinite-dimensional Hamiltonian action, as shown by Donaldson \cite{donaldson} and Fujiki \cite{fujiki}.
The problem of finding {\em constant scalar curvature metrics} (csck) in a fixed K\"ahler class was recently solved through the groundbreaking work of Chen and Cheng \cite{chen-cheng-1, chen-cheng-2}.
Recall that csck metrics are special cases of {\em extremal metrics} \cite{calabi}, namely, K\"ahler metrics that extremize the so-called Calabi functional.
For a thorough introduction to the subject, we refer to \cite{gauduchon-book, szekelyhidi, guedj-zeriahi}.

\section{Hermitian metrics on holomorphic manifolds}
In the Hermitian setting, there is a first issue to address.
Let $X$ be a holomorphic manifold of complex dimension $n$, whose complex structure will be denoted by $J$, and endowed with a Hermitian metric $g$, with associated $(1,1)$-form $\omega=g(J\_,\_)$.
The Levi-Civita connection does not preserve the Hermitian structure, unless it is K\"ahler. This follows from the identities $3d\omega(V,W,Z) = g((D_VJ)W,Z) + g((D_WJ)Z,V) + g((D_ZJ)V,W)$ and $2g((D_VJ)W,Z) = 3d\omega(V,W,Z) - 3d\omega(V,JW,JZ)$.
Then, if we want to take into account the complex structure, we must consider linear connections $\nabla$ that preserve both the tensor of the Hermitian structure, namely $\nabla g=0$ and $\nabla J=0$, even if this comes at the cost of introducing torsion terms $T$.

\subsection{Gauduchon connections}

The space of Hermitian connections on $X$ is an affine space modeled on $\Omega^{1,1}(X,TX)$. More precisely, any Hermitian connection $\nabla$ is determined by $g(\nabla_VW,Z)=g(D_VW,Z)+A$, where the potential $A$ is a vector valued $2$-form, view as a trilinear form that is skew-symmetric in the last two arguments. The potential completely determines the torsion, as detailed in \cite[page 265]{gauduchon-BUMI}; see also detailed computations in \cite{angella-pediconi-RivPr}.
The space $\Omega^2(X,TX)$ of vector-valued $2$-forms, where the torsion lives, admits a natural decomposition $\Omega^2(X,TX) = \Omega^1(X) \oplus (\Omega^2(X,TX))^0 \oplus \Omega^3(X)$.
The corresponding components of an element $B\in\Omega^2(X,TX)$ are denoted as $B = \widetilde{\mathrm{tr}\,B}+B^0+\mathfrak bB$.
Here:
the space $\Omega^3(X)$ of $3$-forms is embedded as a subspace of $\Omega^2(X,TX)$ via $B(V,W,Z)=g(V, B(W,Z))$;
the projection $\mathfrak b \colon \Omega^2(X,TX) \to \Omega^3(X)$ is called the Bianchi projector;
the operator $\mathrm{tr}\, B\colon \Omega^2(X,TX)\to \Omega^1(X)$ is defined as the trace of the bilinear form $(V,W)\mapsto g(V,B(W,\_))$;
the space of $1$-forms $\Omega^1(X)$ is embedded as a subspace of $\Omega^2(X,TX)$ via $\tilde\alpha(V,W,Z)=\sfrac{1}{(2n-1)} (\alpha(Z)g(V,W)-\alpha(W)g(Z,V))$ where $\alpha\in\Omega^1(X)$;
the space $(\Omega^2(X,TX))^0$ consists of trace-free vector-valued $2$-forms $B^0$ satisfying the Bianchi identity $\mathfrak bB^0=0$.
With respect to the complex structure $J$, the space $\Omega^2(X,TX)$ admits the decomposition $\Omega^2(X,TX) = \Omega^{2,0}(X,TX) \oplus \Omega^{1,1}(X,TX) \oplus \Omega^{0,2}(X,TX)$, where: $\Omega^{2,0}(X,TX)$ consists of vector-valued $2$-forms $B$ such that $B(JV,W)=J(B(V,W))$; the space $\Omega^{1,1}(X,TX)$ consists of $B$ such that $B(JV,JW)=B(V,W)$; and $\Omega^{2,0}(X,TX)$ consists of $B$ such that $B(JV,W)=-J(B(V,W))$.
In particular, we consider the orthogonal splitting of $\Omega^{1,1}(X,TX)$ as $\Omega^{1,1}(X,TX)=\Omega^{1,1}_s(X,TX)\oplus\Omega^{1,1}_a(X,TX)$, where $\Omega^{1,1}_s(X,TX)$ denotes the subspace of elements satisfying the Bianchi identity, and $\Omega^{1,1}_a(X,TX)$ is its orthogonal complement, which is isomorphic to $\Omega^3(X)\cap (\Omega^{2,1}(X)\oplus\Omega^{1,2}(X))$.
It turns out, as shown in \cite[Proposition 2]{gauduchon-BUMI}, that the component $T^{0,2}$ does not depend on the connection and is equal to the Nijenhuis tensor, which therefore vanishes. Moreover, the term $\mathfrak b(T^{2,0}-T^{1,1}_a)$ is also independent of the connection and equals the $(2,1)+(1,2)$ component of $\sfrac{1}{3}d^c\omega$, denoted by $\sfrac{1}{3}(d^c\omega)^+$. As a result, the connection is completely determined by $T^{1,1}_s$ and the $(2,0)+(1,1)$ component of $\mathfrak b T$, denoted by $(\mathfrak b T)^+$.
In his celebrated paper \cite{gauduchon-BUMI}, Gauduchon introduced a {\em canonical} one-parameter family $\{\nabla^t\}_{t\in\mathds R}$ of Hermitian connections, now called {\em Gauduchon connections}, by setting $T^{1,1}_s=0$ and $(\mathfrak b T)^+=\sfrac{(2t-1)}{3} (d^c\omega)^+$. More explicitly, they are given by:
$$ g(\nabla^t_V W, Z) = g(D_VW,Z) - \frac{1-t}{4} d\omega(JV,JW,JZ) + \frac{1+t}{4} d\omega(JV,W,Z) . $$
Of course, they all reduce to the Levi-Civita connection if the metric is K\"ahler.
Among them, two distinguished connections are particular noteworthy. For $t=1$, one obtains the {\em Chern connection}, which is the unique Hermitian connection whose component $\nabla^{0,1}$ corresponds to the Cauchy-Riemann operator $\overline\partial_{TX}$ of the holomorphic tangent bundle. For $t=-1$, one obtains the {\em Bismut connection} \cite{bismut}, which is the unique Hermitian connection with totally skew-symmetric torsion. We also notice the following special cases. The connection for $t=0$ is the orthogonal projection of the Levi-Civita connection onto the affine space of Hermitian connections; it is also known as the first canonical connection. The connection for $t=\sfrac{1}{2}$ is referred as conformal by Libermann. The connection for $t=\sfrac{1}{3}$ is characterized by having minimal torsion among Hermitian connections.
Note that, more generally, one can consider the family of metric connections $\nabla^{\varepsilon, \rho}$ defined by
$$ g(\nabla^{\varepsilon,\rho}_V W, Z) = g(D_VW,Z) - \varepsilon d\omega(JV,JW,JZ) + \rho d\omega(JV,W,Z) , $$
for $(\varepsilon,\rho)\in\mathds R^2$, as introduced by \cite{otal-ugarte-villacampa-NuclPhysB}. Among these, when the metric is non-K\"ahler, the only Hermitian connections lie on the line $\varepsilon+\rho=\sfrac{1}{2}$, corresponding to the Gauduchon connections, with the parameter relation $t=1-4\varepsilon$.

\subsection{Curvatures}

Associated with any of these connections, there is a {\em curvature operator}
$$ R \in \Omega^2(X,\mathrm{End}(TX)), \qquad R(V,W):=[\nabla_V,\nabla_W]-\nabla_{[V,W]},$$
(note the different notation used in the literature compared to the Riemannian case,) which we can also be viewed as $R\in\Omega^2(X)\otimes\Omega^2(X)$ by setting $R(V,W,Z,Y):=g(R(V,W)Z,Y)$.
In the K\"ahler case, all Gauduchon connections reduce to the Levi-Civita connection, which also coincides with the Chern connection. As a consequence, the additional symmetries $R\in S^2\Omega^2(X)$ and $R \in \Omega^{1,1}(X,\mathrm{End}(T^{1,0}X))$ hold. The question of when these symmetries still hold in the non-K\"ahler case was first studied by Gray \cite{gray}, and Yang and Zheng \cite{yang-zheng-2}, see also \cite{angella-otal-ugarte-villacampa} for further details.
The missing symmetries give rise to several notions of {\em Ricci curvature}, arising from different traces:
$$ \mathrm{Ric}^{(1)}(V,W)=\mathrm{tr}\,R(V,W)\in\Omega^2X,$$
$$\mathrm{Ric}^{(2)}(Z,Y)=\mathrm{tr}_gR(\_,\_,Z,Y)\in\Omega^2(X)\subset\mathrm{End}(TX),$$
both yielding the same {\em scalar curvature}
$$ s:=\mathrm{tr}_g\mathrm{Ric}^{(1)}=\mathrm{tr}\,\mathrm{Ric}^{(2)} \in \mathcal C^\infty(S,\mathds R).$$
There is also a third way to take a trace, but we have not found any clear geometric motivation for it.
Note that, for the Chern connection,  the first Ricci curvature is a closed real $(1, 1)$-form that represents the first Bott-Chern class $c_1^{BC}(X)\in H^{1,1}_{BC}(X,\mathds R)$, mapped to the usual Chern class of the anti-canonical bundle $K_X^{-1}$.

\subsection{Special metrics}
Let us look at $D\omega$, the covariant derivative of the $(1,1)$-form $\omega$ with respect to the Levi-Civita connection. It can be viewed as an element of $\Omega^2(X,TX)$ via $g((D\omega)(V,W),Z)=(D_V\omega)(W,Z)$. See \cite[Proposition 1]{gauduchon-BUMI} for the determination of each component in the complex bi-degree decomposition.
With respect to the irreducible representations of $U(n)$, the tensor $D\omega$ decomposes as the sum of four components, one of which vanishing in the integrable case. The vanishing of some subsets of the set of four irreducible components led to the celebrated Gray and Hervella classification \cite{gray-hervella} into $2^4=16$ classes of almost-Hermitian metrics.
The K\"ahler case correspond to all the components vanishing, while the general case corresponds to none of them vanishing. Other special classes include: {\em almost-K\"ahler metrics}, namely almost-Hermitian metrics with $d\omega=0$); {\em balanced metrics} \cite{michelsohn}, namely Hermitian metrics satisfying $d\omega^{n-1}=0$, which are particularly interesting as they are invariant under modifications \cite{alessandrini-bassanelli}); {\em locally conformally K\"ahler metrics}, namely Hermitian metrics such that $d\omega=\theta\wedge\omega$ with $d\theta=0$, meaning they locally admit a conformal change to a K\"ahler metric, see, for instance, the recent book \cite{ornea-verbitsky}.

It would be interesting to investigate whether special metrics exist on certain classes of complex manifolds. For instance, on compact complex manifolds satisfying the $\partial\overline\partial$-lemma property, the existence of balanced metrics may follow from solving Hessian equations of the type considered in \cite{szekelyhidi-tosatti-weinkove} or \cite{popovici}.
Note that, when dealing with possibly non-K\"ahler metrics, deformations of $\omega$ of the form $\omega + \sqrt{-1}\partial\overline\partial u$, where $u$ is a smooth real function, do not generally preserve the total volume: see \cite{guedj-lu-1, guedj-lu-2, guedj-lu-3, angella-guedj-lu} for results and applications.

\section{Analytic problems for the Chern connection}

We are now ready to analyze the analogue of the 'constant curvature' problems in the Hermitian setting. For simplicity, we will mainly focus on the Chern connection, although other Hermitian connections are also worth considering.

\subsection{Constant holomorphic sectional curvature}

Let $X$ be a holomorphic manifold endowed with a Hermitian metric $g$. Consider its Chern connection $\nabla$ (or any other Hermitian connection) and the associated Riemann curvature tensor $R$.
A natural object to study is the {\em holomorphic sectional curvature}, defined as
$$ K(V) := K(V,JV) = \frac{g(R(V,JV)V,JV)}{g(V,V)^2} , $$
even though, in general, it no longer encodes the full curvature tensor.

A foundational result in this context is due to Kobayashi: if the holomorphic sectional Chern-curvature $K$ is everywhere negative and bounded away from zero, then the manifold is {\em hyperbolic} in the sense of Kobayashi, see \cite{kobayashi-book}. This result generalizes the Schwarz-Ahlfors-Pick theorem from dimension $1$ to higher dimension. In particular, by the Brody theorem, if $X$ is compact, every holomorphic map $\mathds C \to X$ must be constant.

Complete K\"ahler manifolds with constant holomorphic sectional curvature are referred to as {\em complex space forms}. Analogously to the Riemannian case, their universal covers are classified as follows: the complex projective space $\mathds CP^n$ with the standard Fubini-Study metric, the complex Euclidean space $\mathds C^n$ with the flat metric, and the complex hyperbolic space $\mathds CH^n$ equipped with the Bergman metric; see {\itshape e.g.} \cite[IX.7]{kobayashi-nomizu-2}.

In the non-K\"ahler setting, the first attempt to study compact Hermitian manifolds with constant holomorphic sectional Chern-curvature dates back to \cite{balas}. Examples include compact quotients of complex Lie groups, which are precisely the compact Chern-flat manifolds thanks to \cite{boothby}.
It is conjectured that a compact Hermitian manifold with holomorphic sectional Chern-curvature equal to $c$ must be K\"ahler (hence a complex space form) when $c\neq0$, or Chern-flat when $c=0$.
This conjecture has been confirmed for compact complex surfaces by \cite{balas-gauduchon} (in the case of constant non-positive holomorphic sectional Chern-curvature) and \cite{apostolov-davidov-muskarov} (in the case of pointwise constant holomorphic sectional Chern-curvature); for twistor spaces by \cite{davidov-grantcharov-muskarov}; for locally conformally K\"ahler manifolds with non-positive constant holomorphic sectional Chern-curvature by \cite{chen-chen-nie}; for complex nilmanifolds with left-invariant Hermitian structures by \cite{li-zheng}; for non-balanced Bismut torsion parallel manifolds, including Vaisman manifolds, by \cite{chen-zheng-Chern}; and for Hermitian manifolds whose Bismut connection obeys all the K\"ahler symmetries (so-called Bismut-K\"ahler-like) by \cite{rao-zheng}.
Hermitian metrics with vanishing holomorphic Chern-curvature on compact holomorphic manifolds are studied in \cite{broder-tang}.

The same questions can be posed for Hermitian metrics with constant holomorphic sectional curvature with respect to other connections. The case of the Levi-Civita connection was considered in \cite{sato-sekigawa}, with positive results for compact complex surfaces \cite{apostolov-davidov-muskarov}. The case of the Bismut connection was studied in \cite{chen-zheng-Bismut}.
More generally, an extension to the $2$-parameter family of canonical metric connections, including the Gauduchon connections, was proposed and studied by \cite{chen-nie}, with further developments by \cite{chen-zheng-Canonical}.
Compact Hermitian manifolds admitting a flat Gauduchon connection were investigated by \cite{yang-zheng} and their classification was completed in \cite{lafuente-stanfield}: they are either K\"ahler flat or Chern-flat or Bismut-flat, and the same classification holds when the flatness assumption is replaced by the so-called K\"ahler-like condition, confirming a conjecture by \cite{angella-otal-ugarte-villacampa}.

To conclude this section, we recall another intermediate notion of interest in K\"ahler geometry.
Given two complex planes $\pi=\mathrm{span}\{V,JV\}$ and $\pi'=\mathrm{span}\{W,JW\}$, the {\em holomorphic bisectional curvature} \cite{goldberg-kobayashi} is defined as $H(V,W) := \sfrac{g(R(V,JV)W,JW)}{(g(V,V)g(W,W)-g(V,W)^2)}$.
Since $H(X,X)=K(X)$, the holomorphic bisectional curvature contains more information than the holomorphic sectional curvature.
Moreover, since $H(V,W)=K(V,W)+K(V,JW)$ by the Bianchi identity, it provides less information than the (real) sectional curvature.
The positivity of $H$ is related to the Griffiths positivity \cite{griffiths} of the holomorphic tangent bundle.
In \cite{siu-yau}, Siu and Yau gave an affirmative answer to the Frankel conjecture \cite{frankel}, which states that every compact K\"ahler manifold with positive holomorphic bisectional curvature is biholomorphic to the projective space, serving as an analogue of the sphere theorem in Riemannian geometry. (Using algebraic techniques, Mori \cite{mori} proved the stronger result of the Hartshorne conjecture \cite{hartshorne}. For the case of non-negative holomorphic bisectional curvature, see also \cite{mok}.)
An alternative proof, utilizing the K\"ahler-Ricci flow, was provided by \cite{chen-sun-tian} (see also \cite{chen-AdvMath}).
In the non-K\"ahler setting, Ustinovskiy \cite{ustinovskiy, ustinovskiy-2} proved that, among the Hermitian curvature flow proposed by \cite{streets-tian}, there exists one that preserves Griffiths-non-negativity of the Chern curvature. Combined with the results of \cite{siu-yau, mori}, this extends the Frankel conjecture to the Hermitian context: a compact Hermitian manifold with Griffiths-non-negative Chern curvature, strictly positive at some point, is biholomorphic to the projective space. A similar study for the Griffiths-positivity of the Bismut curvature was undertaken in \cite{barbaro-jgp}.

\subsection{Constant Chern scalar curvature}
Let $X$ be a compact holomorphic manifold of complex dimension $n$. Once fixed a Hermitian metric $\omega$, its conformal class $\{\omega\} := \{ \exp(f)\omega : f\in\mathcal C^\infty(X,\mathds R)\}$ consists entirely of Hermitian metrics.
Motivated by the Yamabe problem in Riemannian geometry, we ask whether the conformal class $\{\omega\}$ contains at least one Hermitian metric with constant scalar curvature with respect to the Chern connection. More precisely, let us study the moduli space
$$ \mathcal{C}h\mathcal{Y}a(X,\{\omega\}) := \sfrac{\left\{ \omega'\in\{\omega\} : s^{Ch}(\omega') \text{ constant} \right\}}{\mathcal{HC}onf(X,\{\omega\}) \times \mathds R^{>0}} $$
of constant Chern scalar curvature metrics, where $\mathcal{HC}onf(X,\{\omega\})$ denotes the group of biholomorphisms preserving the conformal class $\{\omega\}$, and $\mathds R^{>0}$ acts by homotheties.

Note that this problem is essentially different from the classical Yamabe problem. In fact, as shown in \cite[Corollary 4.5]{liu-yang}, on a compact holomorphic manifold, if the (average of the) Chern scalar curvature of a Hermitian metric equals the (average of the) scalar curvature of the corresponding Riemannian metric, then the metric must be K\"ahler. For the same reason, the Chern-Yamabe problem also differs from the Yamabe problem for almost Hermitian manifolds studied by del Rio and Simanca in \cite{delrio-simanca}.

Under conformal changes, the Chern scalar curvature transforms according to the formula:
$$ s^{Ch}(\exp(\sfrac{2f}{n})\omega) = \exp(-\sfrac{2f}{n})\cdot (s^{Ch}(\omega)+\Delta_\omega^{Ch}f) , $$
where $\Delta_\omega^{Ch}$ denotes the Chern Laplacian with respect to $\omega$. It is given by $\Delta_\omega^{Ch}f=-2\sqrt{-1}\mathrm{tr}_\omega\partial\overline\partial f = \Delta_{d,\omega}f+g(df,\theta)$, where $\theta$ is defined by the condition $d\omega^{n-1}=\theta\wedge\omega^{n-1}$ and is called {\em Lee form} or {\em torsion $1$-form}.
Two special cases are of particular interest.
When $\omega$ is {\em balanced} ({\itshape i.e.}, $d\omega^{n-1}=0$, equivalent to $\theta=0$), the operator $\Delta_{\omega}^{Ch}$ coincides with the usual Hodge-de Rham Laplacian $\Delta_{d,\omega}$.
When $\omega$ is {\em Gauduchon} ({\itshape i.e.}, $\partial\overline\partial\omega^{n-1}=0$, equivalent to $d^*_\omega\theta=0$), then the Green identity $\int_X\Delta_\omega^{Ch}f\omega^n=0$ holds.

We recall the following foundational theorem by Gauduchon.
\begin{theorem}[{\cite{gauduchon-cras}}]
On a compact holomorphic manifold, any conformal class of Hermitian metrics contains a unique Gauduchon representative of volume $1$.
\end{theorem}
Thanks to this result, we will fix $\eta$ as the Gauduchon representative of volume $1$ in $\{\omega\}$.

The equation $s^{Ch}(\exp(\sfrac{2f}{n})\omega)=\lambda$, with $\lambda\in\mathds R$, translates into the Liouville-type equation
$$ \Delta_{d,\eta}f + \eta(df,\theta) + s^{Ch}(\eta) = \lambda \cdot \exp(\sfrac{2f}{n}) . $$
Up to homotheties, we choose the normalization $\sfrac{1}{n!}\int_X \exp(\sfrac{2f}{n})\eta^n=1$ for the solution.
The choices for reference and normalization above make the expected constant Chern scalar curvature $\lambda$ completely determined by $X$ and the conformal class $\{\omega\}$:
$$ \lambda = \int_X s^{Ch}(\eta) \eta^n = \frac{1}{(n-1)!} \int_X c_1^{BC}(K_X^{-1}) \wedge [\eta^{n-1}] , $$
where $c_1^{BC}(K_X^{-1}) \in H^{1,1}_{BC}(X,\mathds R)$ maps to the first Chern class of $X$ in $H^2(X,\mathds R)$, and $[\eta^{n-1}]\in H^{n-1,n-1}_A(X,\mathds R)$ is well-defined thanks to $\eta$ being Gauduchon.
We will denote this quantity as the {\em Gauduchon degree} $\Gamma_X(\{\omega\})$.
It is equal to the volume of the divisor associated with any meromorphic section of the anti-canonical line bundle $K_X^{-1}$, measured with respect to the Gauduchon metric $\eta$, see \cite{gauduchon-cras-2, gauduchon-mathann}.
Recall that the {\em Kodaira dimension} is a fundamental holomorphic invariant that measures the growth of the sections of the pluricanonical bundles:
$$ \mathrm{Kod}(X) = \limsup_{m \to +\infty} \frac{\log \dim H^0(X,K_X^{\otimes m})}{\log m} \in \{-\infty, 0, 1, \ldots, n\} .$$
In particular, when $\mathrm{Kod}(X)\geq0$, the Gauduchon degree satisfies $\Gamma_X(\{\omega\})\leq 0$ for any conformal class $\{\omega\}$. More precisely, we have $\Gamma_X(\{\omega\})<0$ for any conformal class, unless $K_X$ is holomorphically torsion, {\itshape i.e.}, there exists an integer $\ell$ such that $K_X^{\otimes \ell}=\mathcal O_X$, in which case $\Gamma_X(\{\omega\})=0$ for any conformal class.

When $\Gamma_X(\{\omega\})=0$, meaning the expected Chern scalar curvature is identically zero, the Chern-Yamabe equation reduces to a linear pde. The Gauduchon condition on $\eta$ guarantees that the kernel of $\Delta_\eta^{Ch})^*$ consists only of constants, and the condition $\Gamma_X(\{\omega\})=0$ assures that $-s^{Ch}$ lies in the orthogonal  complement of $\ker(\Delta_\eta^{Ch})^*)^\perp$, thus making the linear pde solvable.
When $\Gamma_X(\{\omega\})<0$, meaning the expected constant Chern scalar curvature is negative, we can apply conformal techniques and standard elliptic theory, using the continuity method, to prove the existence of a unique solution.
Summarizing, we have the following result, which in particular applies to any conformal class when $\mathrm{Kod}(X)\geq 0$:

\begin{theorem}[{\cite{angella-calamai-spotti}}]
Let $X$ be a compact holomorphic manifol endowed with a conformal class $\{\omega\}$ of Hermitian metrics.
If the Gauduchon degree satisfies $\Gamma_X(\{\omega\})\leq 0$, then the Chern-Yamabe moduli space consists of a single point.
\end{theorem}

When $\Gamma_X(\{\omega\})>0$, the arguments relying on the maximum principle fail, and whether the moduli space $\mathcal Ch\mathcal Ya(X,\{\omega\})$ is non-empty remains an open question.
An important example of this kind is the Hopf surface, which admits a Hermitian metric with positive constant Chern-scalar curvature, see \cite{gauduchon-ivanov}.
Further examples can be constructed using an implicit function theorem argument, starting from a metric whose Chern-scalar curvature is sufficiently small in $\mathcal C^{0,\alpha}$-norm.
We should note that the Chern-Yamabe problem is not variational in general: the Chern-Yamabe equation can be interpreted as the Euler-Lagrange equation of an associated functional if and only if the conformal class contains a balanced representative. 
However, it remains unclear whether this functional is bounded from below.
In the positive case, uniqueness of constant Chern-scalar curvature metrics in a fixed conformal class generally fails. This can be shown by adapting the argument in \cite{delima-piccione-zedda} and employing a version of the Krasnosel’skii Bifurcation Theorem.
We still expect compactness of $\mathcal Ch\mathcal Ya(X,\{\omega\})$, provided the volumes are bounded.
As an additional strategy to attack the existence problem, we introduced the Chern-Yamabe flow. Further attempts and advancements in this direction have been made in \cite{lejmi-maalaoui, calamai-zou, yu}.

Other related problems include the following.
Lejmi and Upmeier \cite{lejmi-upmeier} studied the problem in the (non-integrable) almost-complex setting.
Barbaro \cite{barbaro} addressed the analogous problem for the Bismut connection, and more generally, for connections in the Gauduchon family.
Fusi \cite{fusi} examined the prescribed Chern-scalar curvature problem, analogous to the work of Kazdan and Warner \cite{kazdan-warner} in the Riemannian setting.
In \cite{angella-calamai-pediconi-spotti}, we extended the classical Donaldson-Fujiki interpretation of the scalar curvature as moment map in K\"ahler Geometry to the wider framework of locally conformally K\"ahler Geometry.

\subsection{Chern-Einstein problems}

Regarding the Einstein problem with respect to the Chern connection, two differences immediately appear with the Riemannian case. First, the lack of symmetries in the Chern-Riemann curvature tensor leads to different ways of tracing out the Ricci curvature, resulting in various Chern-Einstein problems. Second, the failure of the first Bianchi identity results in the Einstein factor possibly being non-constant.
We then have three distinct Chern-Einstein problems, each in a weak (when $\lambda$ is a function) or strong (when $\lambda$ is a constant) version:
$$ \mathrm{Ric}^{(i)}(g) = \lambda(x) \cdot g , $$
for $i\in\{1,2,3\}$.
Since the third Chern-Ricci curvature lacks a clear geometric interpretation, we focus our attention on the first and second Chern-Einstein problems.

Let us begin with the {\em first Chern-Einstein problem}.
The case $\lambda=0$ is well understood and corresponds to the so-called {\em non-K\"ahler Calabi-Yau manifolds} \cite{tosatti}, characterized by $c_1^{BC}(X)=0$.

\begin{theorem}[{\cite{tosatti, tosatti-weinkove}}]
Let $X$ be a compact holomorphic manifold endowed with a Hermitian metric $\omega$.
If $c_1(X)=0$, then there exists a first-Chern-Ricci flat metric.
This solution can be expressed either as a conformal transformation of $\omega$, or as a deformation of the form $\omega + \sqrt{-1} \partial\overline\partial \varphi$.
\end{theorem}

In \cite{szekelyhidi-tosatti-weinkove}, the authors prove the Gauduchon’s generalization of the Calabi’s conjecture for compact holomorphic manifolds $X$ with $c_1^{BC}(X)=0$. It follows that that one can always find {\em Gauduchon} first-Chern-Ricci flat metrics, thereby establishing the existence of non-K\"ahler special metrics satisfying both curvature and cohomological conditions.

When the manifold is {\em non-K\"ahler}, meaning that it does not admit any K\"ahler metric, no other interesting cases arise.
This is evident for the {\em strong} problem, as $\mathrm{Ric}^{(1)}(\omega)$ is a closed form. For the {\em weak} problem, conformal arguments yield the following result:

\begin{theorem}[{\cite{angella-calamai-spotti-2}}]
Let $X$ be a compact holomorphic manifold endowed with a Hermitian metric $\omega$ satisfying the first-Chern-Einstein problem with non-identically-zero Einstein factor.
Then $\omega$ is conformal to a K\"ahler metric in the class $\pm c_1(X)$.
\end{theorem}

For this reason, the {\em second Chern-Einstein problem} seems more promising to understand the geometry of non-K\"ahler manifolds.
We observe that, under conformal changes, the second-Chern-Ricci curvature transforms as $\mathrm{Ric}^{(2)}(\exp(f)\omega) = \mathrm{Ric}^{(2)}(\omega)-(\Delta^{Ch}_\omega f)\omega$, in contrast to the first-Chern-Ricci curvature, which transforms as $\mathrm{Ric}^{(1)}(\exp(f)\omega) = \mathrm{Ric}^{(1)}(\omega)-n\sqrt{-1}\partial\overline\partial f$.
In particular, this makes the second-Chern-Einstein problem depends only on the conformal class \cite{gauduchon-cras1980}.
As a consequence, conformal methods allow us to assume that the Einstein factor has a definite sign \cite{gauduchon-cras1981}, equal to the sign of the Gauduchon degre and determined according to $K_X$ and $K_X^{-1}$ being pseudo-effective or unitary flat \cite{teleman, yang}. Furthermore, if the Chern-Yamabe problem can be solved affirmatively, then it would even be possible to make the Einstein factor constant.
Obstructions {\itshape à la} Bochner arise from \cite{gauduchon-bsmf, liu-yang-2}: compact manifolds admitting a positive second-Chern-Einstein metric do not have non-trivial holomorphic $p$-forms, while those admitting negative second-Chern-Einstein manifolds do not have non-trivial holomorphic $p$-vector-fields, for $p \geq 1$.
Moreover, second-Chern-Ricci metrics $g$ are weakly-$g$-Hermite-Einstein, see {\itshape e.g.} \cite{lubke-teleman}. This provides additional obstructions: the Bogomolov-L\"ubke inequality holds \cite{bogomolov, lubke}, and the Kobayashi-Hitchin correspondence assures that the holomorphic tangent bundle is $g$-semi-stable \cite{kobayashi, lubke-2}.
The fundamental example in this context is provided by Hopf surfaces, which admit second-Chern-Einstein metrics with positive scalar curvature. Indeed, they are the only non-K\"ahler compact complex surface that admit second-Chern-Einstein metrics \cite{gauduchon-ivanov}.
Other compact examples include holomorphically-parallelizable manifolds \cite{boothby} and the homogeneous non-K\"ahler $C$-spaces studied by Podestà \cite{podesta}. See also \cite{angella-pediconi} for some non-compact simply-connected examples of cohomogeneity-one.
Currently, no compact non-K\"ahler examples admitting second-Chern-Einstein metric with negative scalar curvature are known \cite[Remark 1]{angella-calamai-spotti-2}.

For the analogous problem in the context of almost-K\"ahler structures within the non-integrable setting, we refer to {\itshape e.g.} \cite{apostolov-draghici, dellavedova, barbaro-lejmi}.
For the corresponding problem involving the Bismut connection, we refer to \cite{barbaro, barbaro-2}, which focus on {\em Calabi-Yau with torsion} metrics, meaning Bismut-Ricci-flat metrics. This condition is of particular interest because coupling the equation $(\mathrm{Ric}^B(\omega))^{1,1} = \lambda\omega$, where $\lambda$ is a real function, with the pluriclosed condition $\partial\overline\partial \omega=0$ leads to either K\"ahler-Einstein metrics (when $\lambda\neq0$) or to the so-called {\em Bismut Hermitian Einstein metrics} \cite{garciafernandez-streets} (meaning pluriclosed and Calabi-Yau with torsion), which are stationary points of the pluriclosed flow \cite{streets-tian-IMRN}.
We point out that the pluriclosed flow is motivated by, and conjecturally related to, the geometrization of compact complex surfaces \cite{streets}.
Moreover, Bismut Hermitian-Einstein metrics are closely connected to generalized geometry, see \cite{garciafernandez-streets}.
In particular, the pluriclosed flow preserves generalized K\"ahler geometry in the sense of  Hitchin and Gualtieri \cite{streets-tian-NuclPhysB}, and there are interesting results concerning its global existence and convergence on compact Bismut-flat manifolds \cite{garciafernandez-jordan-streets}. Recently, the authors of \cite{apostolov-barbaro-lee-streets} proved that a compact Hermitian threefold endowed with a pluriclosed metric with vanishing Bismut-Ricci form is either K\"ahler Calabi-Yau or Bismut-flat.

Another interesting generalization in the non-K\"ahler setting is provided by the {\em Hull-Strominger system} \cite{strominger, hull}, see also the recent survey \cite{garciafernandez-survey}. Let $X$ be a complex manifold of complex dimension $n=3$, endowed with a never vanishing holomorphic global section $\Omega$ trivializing the canonical bundle $K_X$. Denote by $M$ the underlying $6$-dimensional smooth manifold endowed with the (integrable) almost-complex structure $J$ and consider the smooth complex vector bundle $(TM, J)$. Let $\mathcal E$ be a holomorphic bundle on $X$. The Hull-Strominger system, with the heterotic equation of motion \cite{ivanov}, aims at finding a Hermitian metric $g$ on $X$, a Hermitian metric $h$ on the bundle $\mathcal E$, and an integrable Dolbeault operator $\overline\partial_{TM}$ on $(TM,J)$ satisfying the Hermite-Einstein equations
$$ \Lambda_\omega F_h=0, \qquad \Lambda_\omega R=0, $$
where $F_h$ denotes the curvature of $(\mathcal E, h)$, respectively $R$ the curvature of the Chern connection on $(TM,J,\overline\partial_{TM})$, and the dilatino equation $d^*\omega - d^c\log\|\Omega\|_\omega = 0$, which is equivalent to the conformally balanced equation \cite{gauntlett-martelli-waldram, li-yau}
$$ d(\|\Omega\|_\omega\omega^{n-1})=0, $$
and the Green-Schwarz anomaly cancellation condition, also known as Bianchi identity,
$$ dd^c\omega-\alpha(\mathrm{tr}\, R\wedge R-\mathrm{tr}\,F_h\wedge F_h)=0,$$
for some parameter $\alpha\in\mathds R$, which is usually taken positive for physical reasons, where the first term denotes the differential of the torsion $-d^c\omega$ of the Bismut connection of $(M,J,g)$.
Alternatively, one can fix a Hermitian metric $h$ on $E$ and rephrase the Hull-Strominger system as coupling a pair of Hermite-Yang-Mills connections, $A$ on $(E,h)$ and $\nabla$ on $(TM,J,\omega)$, meaning
$$ \Lambda_\omega F_A=0, \quad F_A^{0,2}=0 , \qquad \Lambda_\omega R_\nabla=0, \quad R_\nabla^{0,2}=0, $$
with a conformally balanced metric $\omega$ by means of the Bianchi identity.
Explicit invariant solutions to the Strominger system (with the heterotic equations of motion) are found in \cite{fernandez-ivanov-ugarte-vassilev, otal-ugarte-villacampa-NuclPhysB}
on compact non-K\"ahler homogeneous spaces, obtained as the quotient by a lattice of maximal rank of a nilpotent Lie group, the semisimple group $\mathrm{SL}(2,\mathds C)$ and a solvable Lie group.

We conclude this section by emphasizing the role of Hopf manifolds in non-K\"ahler geometry, serving as a natural counterpart to projective space in the realm of K\"ahler geometry. In \cite{angella-miceli-placini}, we investigated the problem of approximating compact regular (respectively, quasi-regular) Vaisman metrics by metrics induced by immersions (respectively, embeddings) into Hopf manifolds. This provides a non-K\"ahler analogue of the Tian approximation theorem for projective manifolds \cite{tian-2, ruan, zelditch} in a specific non-Kähler context.

\subsection{Hermitian flows}
As already remarked, the Ricci flow proved to be a very powerful tool in differential geometry. When the starting metric is K\"ahler, the evolving metrics remain K\"ahler: in this case, the flow is called the K\"ahler-Ricci flow. On the other hand, if the starting metric is just Hermitian, the flow {\em does not}, in general, preserve Hermitianity, unless the curvature of the Levi-Civita connection of the initial metric satisfies all the K\"ahler symmetries \cite{angella-otal-ugarte-villacampa}. For this reason, several geometric flows have been proposed for holomorphic manifolds to assure that Hermitianity is preserved. Among these are the {\em Chern-Ricci flow} \cite{gill, tosatti-weinkove-CRF}, the {\em Hermitian curvature flows} \cite{streets-tian} including the {\em pluriclosed flow} and the flow by Ustinovskiy \cite{ustinovskiy}, the {\em anomaly flows} \cite{phong-picard-zhang-MathZ, phong-picard-zhang-CAG, picard-LNM}.

Let us focus on the Chern-Ricci flow, originally introduced by Gill \cite{gill} to study compact non-K\"ahler Calabi-Yau manifolds and subsequently investigated in depth by Tosatti and Weinkove \cite{tosatti-weinkove-CRF, tosatti-winkove-surfacesCRF, tosatti-winkove-surveyCRF}. This flow evolves a Hermitian metric $\omega_0$ on a (compact) holomorphic manifold $X$ according to the equation
$$ \frac{\partial \omega(t)}{\partial t} = -\mathrm{Ric}^{Ch}(\omega(t)), \quad \text{ with } \omega(0)=\omega_0 , $$
where $\mathrm{Ric}^{Ch}\stackrel{\text{loc}}{=}-\sqrt{-1}\partial\overline\partial\log\det\omega$ denotes the first Ricci curvature with respect to the Chern connection.
One might expect the behaviour of the Chern-Ricci flow to reflect the underlying holomorphic geometry of $X$, making it particularly useful in completing the classification of compact non-K\"ahler complex surfaces.

We recall that the Enriques-Kodaira-Siu classification (see, for instance, the comprehensive book \cite{barth-hulek-peters-vandeven}) classifies compact complex surfaces according to their Kodaira dimension, which can take the values $-\infty$, $0$, $1$, or $2$.
By the Castelnuovo contraction theorem, we can always contract the rational curves of self-intersection $-1$. Therefore, in the following, we will assume the surface to be {\em minimal}.
Surfaces $S$ of general type, that is, those with $\mathrm{Kod}\,S=2$, are always algebraic.
Focusing on non-K\"ahler surfaces, which are topologically characterized by the first Betti-number being odd \cite{lamari, buchdahl}, only the following cases remain.
For $\mathrm{Kod}\,S = 1$, we have minimal non-K\"ahler properly elliptic surfaces, which are described, up to finite covers, as elliptic bundles, see {\itshape e.g.} \cite{brinzanescu, brinzanescu-book}.
For $\mathrm{Kod}\,S=0$, we have Kodaira surfaces, which have $c_1^{BC}(S)=0$, meaning they are non-K\"ahler Calabi-Yau.
Lastly, the so-called {\em class VII} consists of compact complex surfaces that have $\mathrm{Kod}\,S=-\infty$ and $b_1(S)=1$.
Surfaces $S$ in class VII with $b_2(S)=0$ are either Hopf surfaces or Inoue-Bombieri surfaces \cite{inoue, bombieri}, as proven in \cite{bogomolov, li-yau-zheng, teleman-IJM}. The first examples of surfaces $S$ in class VII with $b_2(S)>0$ were constructed in \cite{inoue-1, inoue-2} and later generalized by Kato \cite{kato}, through an iterative process of blowing up the standard ball in $\mathds C^2$ followed by holomorphic surgery.
These surfaces are characterized by the existence of a global spherical shell, meaning an open subset $U \subset S$ biholomorphic to a neighbourhood of $S^3$ in $\mathds C^2$ such that $S\setminus U$ is connected. Moreover, they are degenerations of blown-up primary Hopf surfaces \cite{nakamura-Sugaku}, and always admit singular holomorphic foliations \cite{dloussky-oeljeklaus}. Every Kato surface $S$ has exactly $b_2(S)$ rational curves, and conversely, every minimal compact complex surface in class VII with exactly $b_2(S)>0$ rational curves is a Kato surface \cite{dloussky-oeljeklaus-toma}. It is expected that no other surfaces exist in class VII beyond these, a conjecture known as the Global Spherical Shell conjecture.
The ``optimistic conjecture'' in \cite[Section 4.6]{fang-tosatti-weinkove-zheng} aims at recovering a global spherical shell from the limiting behavior of the Chern-Ricci flow.

For non-K\"ahler properly elliptic surfaces, the Chern-Ricci flow starting at any Gauduchon metric exists for all non-negative time and exhibits the following collapsing behaviour: the normalized metrics $\sfrac{\omega(t)}{t}$ converge in the sense of Gromov-Hausdorff to to a Riemann surface endowed with the distance function induced by an orbifold K\"ahler-Einstein metric, see \cite{tosatti-weinkove-MathAnn}.
When $\mathrm{Kod}\,S=0$, the behaviour of the Chern-Ricci flow has been studied by Gill \cite{gill}. It was shown that the Chern-Ricci flow exists for all time and converges to a Chern-Ricci flat metric.
For classical Hopf surfaces $\sfrac{\mathds C^2 \setminus 0}{\mathds Z}$, where $\mathds Z$ is generated by $\lambda \cdot\mathrm{id}$ with $\lambda\in\mathds C$ and $0<|\lambda|<1$, the solution starting from the standard metric can be explicitly written \cite{tosatti-winkove-surfacesCRF}. The flow collapses in finite time and smoothly converges to a non-negative $(1,1)$-form, whose kernel defines a smooth distribution, with iterated brackets generating the tangent space of the Sasaki $S^3$. In this sense, the limiting behaviour of the Chern-Ricci flow clearly detects both the differential and holomorphic geometries of the Hopf surface.
See also \cite{edwards} for further results on the Chern-Ricci flow for Hopf surfaces, supporting the conjecture that the Gromov-Hausdorff limit is isometric to a round $S^1$.
The Chern-Ricci flow on Inoue-Bombieri surfaces is well understood after \cite{tosatti-winkove-surfacesCRF, lauret, lauret-valencia, fang-tosatti-weinkove-zheng, angella-tosatti}.

We briefly recall the results on the convergence of the normalized Chern-Ricci flow $\sfrac{\partial}{\partial t}\omega(t)=-\mathrm{Ric}^{Ch}(\omega(t))-\omega(t)$ starting from a Gauduchon metric $\omega(0)=\omega_0$ on an Inoue-Bombieri surface $S$.
First, recall that Inoue-Bombieri surfaces are classified into three types: $S_M$, $S^+$ and $S^-$, each depending on some parameters. Inoue-Bombieri surfaces in class $S_M$ have a differentiable torus-bundle structure over $S^1$ and a holomorphic parabolic foliation induced by a non-negative $(1,1)$-form $\alpha\in -c_1^{BC}(S)$, whose leaves are dense in the fibres of the torus-bundle. They also admit a Hermitian metric, known as the Tricerri metric \cite{tricerri-Torino}, which is flat along the leaves of the foliation and lifts to a K\"ahler metric on the universal covering $\mathds C \times \mathds H$. In \cite{tosatti-winkove-surfacesCRF}, the explicit solution starting from the Tricerri metric $\omega_T$ is computed and analyzed. It is shown that the flow uniformly converges to $\alpha$ as $t\to+\infty$, and the corresponding metric spaces converge in Gromov-Hausdorff to $S^1$. This result is extended when starting from any locally homogeneous metric \cite{lauret-TAMS, lauret-valencia}, and, more generally, from a metric of type $\omega_{SFL}+\sqrt{-1}\partial\overline\partial u$ with $\omega_{SFL}$ strongly flat along the leaves \cite{fang-tosatti-weinkove-zheng}. The latter condition means that $\omega_{SFL}$ restricted to any leaf $\mathcal L\simeq \mathds C$ of the foliation is a K\"ahler flat metric and moreover $p^*\omega_{SFL}=c\cdot\sqrt{-1}dz\wedge d\bar z$, with $c=c(\Im w)$ smooth, where $p\colon \mathds C \times \mathds H \to S$ is the universal covering and $(z,w)\in\mathds C \times\mathds H$ are coordinates, with $\Im w>0$. This is equivalent to the function $u$ solving the degenerate elliptic equation
$$ \frac{\sqrt{-1}\partial\overline\partial u \wedge \alpha}{\omega^2_T} = -\frac{\omega\wedge\alpha}{\omega_T^2} + \frac{\int_S\omega\wedge\alpha}{\int_S\omega_T^2} , $$
which involves the leafwise Laplacian $\Delta_{\mathcal D}=\sfrac{\sqrt{-1}\partial\overline\partial \_\wedge\alpha}{\omega_T^2}$.
Finally, in \cite{angella-tosatti}, we proved that any Gauduchon metric admits a $\sqrt{-1}\partial\overline\partial$-deformation that is strongly flat along the leaves. This is achieved by solving the equation $\Delta_{\mathcal D}u=g$ on $S$.
Exploiting Fourier expansion along the torus fibres (along with the irrationality arising from the parameters in the construction of the surface), we obtain a distributional solution. The Liouville theorem on rational approximation of irrational algebraic
numbers assures that this solution is smooth on each torus. By continuous dependence and uniqueness, we get that the solution is smooth on the torus bundle and then descends to $S$. The same argument applies to surfaces of type $S^+$, which also have a bundle structure with fibres given by quotients of the three-dimensional Heisenberg group and a holomorphic foliation with leaves isomorphic to $\mathds C\setminus\{0\}$, whose leaves are dense in the fibres. In this case, we use partial Fourier expansion on nilmanifolds, a technique that has proven useful in many other contexts, see {\itshape e.g.} \cite{auslander-tolimieri, deninger-singhof, richardson, holt-zhang-MRL, holt-zhang-AdvMath, ruberman-saveliev}, and which beautifully connects complex geometry, pdes, and number theory. Finally, surfaces of type $S^-$ admit an unramified double cover of type $S^+$. Therefore we have the following. (For further results on when the convergence is indeed smooth, see \cite[Theorem 1.3]{angella-tosatti}.)
\begin{theorem}[{\cite{fang-tosatti-weinkove-zheng, angella-tosatti}}]
Let $S$ be an Inoue-Bombieri surface endowed with a Gauduchon metric $\omega$. The solution $\omega(t)$ of the normalized Chern-Ricci flow starting at $\omega$ exists for all time and uniformly converges to a multiple of the degenerate metric $\alpha$. Moreover, the corresponding metric spaces converge to $S^1$ with the round metric in Gromov-Hausdorff sense.
\end{theorem}

\bibliographystyle{alpha}
\bibliography{biblio}

\end{document}